\magnification=1200
\def\defeq{\buildrel \rm def \over =}   

\catcode`\@=11
\font\tenmsb=msbm10
\font\sevenmsb=msbm7
\font\fivemsb=msbm5
\newfam\msbfam
\textfont\msbfam=\tenmsb  \scriptfont\msbfam=\sevenmsb
  \scriptscriptfont\msbfam=\fivemsb
\def\msb@{\hexnumber@\msbfam}
\def\Bbb{\relax\ifmmode\let\next\Bbb@\else
 \def\next{\errmessage{Use \string\Bbb\space only in math mode}}\fi\next}
\def\Bbb@#1{{\Bbb@@{#1}}}
\def\Bbb@@#1{\fam\msbfam#1}

\font\teneuf=eufm10
\font\seveneuf=eufm7
\font\fiveeuf=eufm5
\newfam\euffam
\textfont\euffam=\teneuf \scriptfont\euffam=\seveneuf
  \scriptscriptfont\euffam=\fiveeuf
\def\euf@{\hexnumber@\euffam}
\def\Frak{\relax\ifmmode\let\next\Frak@\else
 \def\next{\errmessage{Use \string\Frak\space only in math mode}}\fi\next}
\def\Frak@#1{{\Frak@@{#1}}}
\def\Frak@@#1{\fam\euffam#1}

\catcode`\@=12
\def\N{{\Bbb N}}

\newcount\thmno \thmno=1
\long\def\theorem#1#2{\edef#1{Theorem~\the\thmno}\noindent
    {\sc Theorem \the\thmno.\enspace}{\it #2}\endgraf
    \penalty55\global\advance\thmno by 1}
\newcount\lemmano \lemmano=1
\long\def\lemma#1#2{\edef#1{Lemma~\the\lemmano}\noindent
    {\sc Lemma \the\lemmano.\enspace}{\it #2}\endgraf
    \penalty55\global\advance\lemmano by 1}
\newcount\propno \propno=1
\long\def\proposition#1#2{\edef#1{Proposition~\the\propno}\noindent
    {\sc Proposition \the\propno.\enspace}{\it #2}\endgraf
    \penalty55\global\advance\propno by 1}
\newcount\corolno \corolno=1
\long\def\corollary#1#2{\edef#1{Corollary~\the\corolno}\noindent
    {\sc Corollary \the\corolno.\enspace}{\it #2}\endgraf
    \penalty55\global\advance\corolno by 1}


\newcount\conjno \conjno=1
\long\def\conjecture#1#2{\edef#1{Conjecture~\the\conjno}\noindent
    {\sc Conjecture \the\conjno.\enspace}{\it #2}\endgraf
    \global\advance\conjno by 1}

%

\newcount\eqnno \eqnno=0
\def\eqn#1{\global\edef#1{(\the\secno.\the\eqnno)}#1\global
   \advance\eqnno by 1}
\newcount\secno \secno=0
\def\sec#1\par{\eqnno=1\global\advance\secno by 1\bigskip
    \noindent\centerline{\sc\the\secno. #1}\par\nobreak\noindent}


\def\proof{\noindent{\it Proof.\enspace}}

%

\def\qed{~\vrule width4pt height6pt depth1pt}

\parskip=12pt plus4pt minus4pt

\font\sc=cmcsc10
\def\itc#1{{\it #1\/}}
\def\sgn{{\rm sgn}\,}
\def\tpof{$({\bf 3+1})$-free}
\def\suml{\sum_{\lambda:\ell(\lambda)=\ell}}

\centerline{\bf A Note on a Combinatorial Interpretation of the}
\centerline{\bf e-Coefficients of the Chromatic Symmetric Function}
\bigskip
\centerline{Timothy Y. Chow}
\medskip
\centerline{Dept.\ of Mathematics, Univ.\ of Michigan,
Ann Arbor, MI 48109, U.S.A.}
\centerline{email: \tt tchow@umich.edu}

\bigskip
{\narrower\noindent
{\bf Abstract.}
Stanley has studied a symmetric function generalization~$X_G$ of the
chromatic polynomial of a graph~$G$.
The innocent-looking Stanley-Stembridge Poset Chain Conjecture states
that the expansion of~$X_G$ in terms of elementary symmetric functions
has nonnegative coefficients if $G$ is a clawfree incomparability graph.
Here we give a combinatorial interpretation of these coefficients by
combining Gasharov's work on the
conjecture with E\~gecio\~glu and Remmel's combinatorial interpretation
of the inverse Kostka matrix.
This gives a new proof of a partial nonnegativity result of Stanley.
As an interesting byproduct we derive a previously
unnoticed result relating acyclic orientations to $P$-tableaux.

}

\bigskip

\sec Introduction

The main ideas in this note are simple
but require an inordinate number of definitions to state.
In this section we skip most of these definitions
so as not to obscure the exposition with a mass of technicalities.
The missing definitions are given in the next section.

Let $G$ be a finite simple undirected graph and let
$X_G$ be its chromatic symmetric function.
Expand $X_G$ in terms of elementary symmetric functions~$e_\lambda$
and call the coefficients $a^G_\lambda$:
$$X_G = \sum_\lambda a^G_\lambda e_\lambda.$$
One of the outstanding conjectures about~$X_G$
(and the motivation for this note)
is the Stanley-Stembridge \itc{Poset Chain Conjecture}~[6]:
if $G$ is a clawfree incomparability graph,
then $a^G_\lambda\ge 0$ for all~$\lambda$.

It is natural to attack this conjecture
by looking for a combinatorial interpretation of~$a^G_\lambda$.
We can obtain such an interpretation as follows.
Observe first that Gasharov~[2] tells us that if $G$ is a
clawfree incomparability graph,
then the coefficients of the Schur-function expansion
of~$X_G$ have a combinatorial interpretation.
Next, observe that to convert from
the Schur-function expansion of~$X_G$
to the $e$-expansion of~$X_G$,
we need to introduce an inverse Kostka matrix.
But E\~gecio\~glu and Remmel~[1]
have a combinatorial interpretation
of the inverse Kostka matrix.
Combining these two combinatorial interpretations therefore gives us
a combinatorial interpretation of the coefficients~$a^G_\lambda$.

This simple observation
does not in itself prove the Poset Chain Conjecture,
because E\~gecio\~glu and Remmel's combinatorial interpretation
(and therefore our combinatorial interpretation of~$a^G_\lambda$)
involves a \itc{signed} sum over combinatorial objects.
However, it does open up a new line of attack on the Poset Chain Conjecture:
we can try to prove the nonnegativity of~$a^G_\lambda$ by looking
for sign-reversing involutions.
This is illustrated below by a new proof of the fact
(first shown by Stanley [5, Theorem~3.3])
that if $G$ is a clawfree incomparability graph,
then for all~$\ell$,
$$\suml a^G_\lambda$$
is the number of acyclic orientations of~$G$ with $\ell$~sinks.
This new proof is not significantly shorter than Stanley's,
but in addition to using completely different methods,
it seems to require significantly less ingenuity;
the sign-reversing involution is very simple and one of
the first things one might try.
This provides hope that more sophisticated involutions
will produce correspondingly stronger results.

An interesting byproduct of the proof is Lemma~1 below,
which describes a connection, apparently not previously noticed,
between acyclic orientations and $P$-tableaux.

\sec Background

We now provide the necessary technical background.
Some familiarity with the basics of symmetric functions and
partitions is assumed; see [4] or [3, Chapter~I].

For the expert, we remark
that there are two points where we diverge
slightly from the literature:
we use English style for our Ferrers diagrams
while E\~gecio\~glu and Remmel use French style,
and by ``$P$-tableau'' we mean the \itc{transpose}
of what Gasharov calls a $P$-tableau.

Let $G$ be a finite simple undirected graph with vertex set
$V = \{v_1, v_2, \ldots, v_d\}$.
A \itc{proper coloring} of~$G$ is a map
$\kappa: V \to \N$ such that $\kappa(v_i)\ne \kappa(v_j)$
whenever $v_i$ and~$v_j$ are adjacent.
Let $\{x_n \mid n \in \N\}$ be a countably infinite family of
independent indeterminates.
Following~[5], define the \itc{chromatic symmetric function~$X_G$}
of~$G$ to be the formal power series
$$X_G \defeq \sum_\kappa
    x_{\kappa(v_1)} x_{\kappa(v_2)} \cdots x_{\kappa(v_d)},$$
where the sum is over all proper colorings of~$G$.
It is easy to see that $X_G$ is a symmetric function,
so it can be written as a (finite) linear combination of
elementary symmetric functions~$e_\lambda$:
$$X_G = \sum_{\lambda\vdash d} a^G_\lambda e_\lambda.$$

As we mentioned in the introduction,
one of the main open problems in this area
is the \itc{Poset Chain Conjecture} of Stanley and Stembridge~[6].
This states that if $G$ is a clawfree incomparability graph, then
$G$ is $e$-positive, i.e., $a^G_\lambda\ge 0$ for all~$\lambda$.
Recall that an \itc{incomparability graph} is a graph obtained from a
finite poset by letting the vertex set of the graph be the vertex
set of the poset and connecting two vertices with an edge if and
only if the vertices are incomparable elements in the poset.
\itc{Clawfree} just means that the graph does not contain
the complete bipartite graph~$K_{1,3}$ as an induced subgraph.
Clawfree incomparability graphs are also referred to as
``incomparability graphs of \tpof\ posets.''
A poset is \itc{\tpof} if it does not contain an induced subposet
isomorphic to
a disjoint union of a three-element chain and a one-element chain.
It is clear that an incomparability graph of a poset is clawfree
if and only if the poset is \tpof.

One of the most important partial results towards the Poset
Chain Conjecture is due to Gasharov.
If $P$ is a finite poset, then define a \itc{$P$-tableau}
to be an arrangement of the elements of~$P$ into a Ferrers shape
(English style) such that

\item{1.} each element of~$P$ is used exactly once,

\item{2.} if $x$ appears immediately above~$y$ in a column
then $x\prec y$, and

\item{3.} if $x$ appears immediately to the left of~$y$ in a row
then $x \not\succ y$.

\noindent
(We have chosen to use the transpose of Gasharov's definition
of $P$-tableaux both for convenience in our proof and
because in general column-strict tableaux are more commonly used
in the literature than row-strict tableaux.)
Next, define $f^G_\lambda$ by the equation
$$\omega X_G = \sum_{\lambda\vdash d} f^G_\lambda s_\lambda,$$
i.e., expand~$\omega X_G$ in terms of Schur functions
and let $f^G_\lambda$ be the coefficient of~$s_\lambda$.
(Here $\omega$ is the involution that
sends $s_\lambda$ to~$s_{\lambda'}$.)
We can now state Gasharov's result~[2].

\proposition\gasharov
{If $P$ is a \tpof\ poset and
$G$ is its incomparability graph, then $f^G_\lambda$ is
the number of $P$-tableaux of shape~$\lambda$.}

The last piece of background is
E\~gecio\~glu and Remmel's
combinatorial interpretation of the inverse Kostka matrix.
A \itc{special rim hook tabloid~$T$} of shape~$\mu$
and type $\lambda = (\lambda_1, \ldots, \lambda_\ell)$
is a filling of the Ferrers diagram
of~$\mu$ repeatedly with rim hooks (a.k.a.\ \itc{skew hooks}
or \itc{border strips} or \itc{ribbons}) of sizes
$\{\lambda_1, \ldots, \lambda_\ell\}$ such that each
rim hook is \itc{special,} by which we mean that it
contains at least one cell in the first column.
In other words, to create a special rim hook tabloid,
take any rim hook that contains at least one cell in the
first column and that leaves a legal Ferrers diagram when
removed; then remove this rim hook and repeat the process
iteratively on the residue.
Note that only the \itc{sizes} of the rim hooks matter
and not their \itc{order,} in contrast to the usual
notion of a rim hook tableau.  In other words,
we may if we wish peel off a large rim hook first,
then a small one, then a large one, and so on,
so long as in the end we have the right number of rim hooks of each size.
For this reason, special rim hook tabloids are typically drawn not by
putting numbers in the boxes of the Ferrers diagram
but by connecting the boxes in question with a continuous zigzag line.

The \itc{sign}~$\sgn T$ of a special rim hook tabloid~$T$
is defined in the expected way:
the sign of a rim hook is $(-1)^{h-1}$ where $h$ is the height
of the rim hook, and the sign of~$T$ is the product of the signs
of its component rim hooks.

E\~gecio\~glu and Remmel's result is the following~[1].

\proposition\egrem{
The inverse Kostka matrix $K^{-1}_{\lambda,\mu}$ satisfies
$$K^{-1}_{\lambda,\mu} = \sum_T \sgn T,$$
where the summation is over all special rim hook tabloids of type~$\lambda$
and shape~$\mu$.}

\sec A Combinatorial Interpretation of $a^G_\lambda$

Throughout this section, unless otherwise stated,
$P$ will be a \tpof\ poset,
$G$ will be its incomparability graph,
and $X_G$ will be the chromatic symmetric function of~$G$.
The coefficients $a^G_\lambda$ and~$f^G_\lambda$
are defined as above.

It is well known (e.g., [3, \S I.6, Table~1]) that
the change-of-basis matrix between the Schur functions
and the complete homogeneous symmetric functions is
the inverse Kostka matrix, i.e.,
$$s_\mu = \sum_\lambda K^{-1}_{\lambda,\mu} h_\lambda.$$
Therefore
$$\sum_\lambda a^G_\lambda h_\lambda
   = \omega X_G = \sum_\mu f^G_\mu s_\mu
   = \sum_\mu f^G_\mu \sum_\lambda K^{-1}_{\lambda,\mu} h_\lambda
   = \sum_\lambda
       \biggl( \sum_\mu K^{-1}_{\lambda,\mu} f^G_\mu \biggr) h_\lambda,$$
i.e.,
$$a^G_\lambda = \sum_\mu K^{-1}_{\lambda,\mu} f^G_\mu.$$
Now let us invoke \egrem.  We obtain
$$\eqalignno{a^G_\lambda &= \sum_\mu \sum_T (\sgn T) f^G_\mu,&\eqn\comb\cr}$$
where the inner sum is over all special rim hook tabloids of type~$\lambda$
and shape~$\mu$.

The right-hand side of~\comb\ suggests the following definition.
A \itc{special rim hook $P$-tableau~$T$} of type~$\lambda$
and shape~$\mu$ is an ordered pair $(T',T'')$
where $T'$ is a $P$-tableau  of shape~$\mu$
and $T''$ is a special rim hook tabloid of type~$\lambda$
and shape~$\mu$.
The \itc{sign} of~$T$ is just the sign of~$T''$.
Note that because $T$ and~$T'$ have the same shape,
we may visualize a special rim hook $P$-tableau as
an ordinary $P$-tableau equipped with a decomposition
into special rim hooks, i.e., we need not visualize two separate
tableaux.

We can now state our combinatorial interpretation of~$a^G_\lambda$.
Combining \gasharov\ and \comb\ yields

\theorem\mainthm
{If $G$ is the incomparability graph of a \tpof\ poset~$P$, then
$$a^G_\lambda = \sum_T \sgn T,$$
where the sum is over all special rim hook $P$-tableaux~$T$
of type~$\lambda$.}

To illustrate the power of \mainthm, we use it to prove

\proposition\acyclic
{If $G$ is the incomparability graph of a \tpof\ poset~$P$,
then for all~$\ell$,
$$\suml a^G_\lambda$$
is the number of acyclic orientations of~$G$ with
exactly $\ell$ sinks.
}

\noindent
Recall that an \itc{acyclic orientation} of~$G$ is an
assignment of a direction to each edge of~$G$ in such a way
that no directed cycles are formed.
Stanley [5, Theorem~3.3] originally proved \acyclic\
with no restriction on~$G$.
For our proof, we need the following lemma.

\lemma\ptableau{Let $G$ be the incomparability graph of
an arbitrary finite poset~$P$.
Let $\kappa_\ell$ be the number of acyclic orientations of~$G$
with exactly $\ell$~sinks and let $\pi_k$ be the number of
$P$-tableaux whose shape is a hook with $k$ cells
in the first column.  Then
$$\pi_k = \sum_\ell {\ell-1 \choose k-1} \kappa_\ell.$$
}

\proof
Let $T$ be a $P$-tableau whose shape is a hook with $k$ cells
in the first column.  Then $T$ induces an acyclic orientation
of $G$ as follows: if $u$ and~$v$ are connected by an edge in~$G$,
then we make $u$ point towards~$v$ if, in~$T$, the column that $u$ is in
lies to the right of the column that $v$ is in, and we make $v$ point
towards~$u$ if the column that $u$ is in lies to the left of the column 
that $v$ is in.
Note that $u$ and~$v$ cannot be in the same column, for then they would
be comparable in~$P$ and therefore non-adjacent in~$G$.
It is clear that this orientation is acyclic.

Let $\Frak o$ be an acyclic orientation of~$G$ with $\ell$ sinks.
We claim that for all~$k$,
$\Frak o$ is induced by exactly $\ell-1 \choose k-1$
$P$-tableaux whose shape is a hook with $k$ cells in the first column.
This will prove the lemma.

Suppose we are given $\Frak o$ and~$k$.
Think of~$\Frak o$ as a poset, with $u<v$ if and only if there is
a directed path from $v$ to~$u$.
To avoid confusing this partial ordering with the partial
ordering of~$P$, we use $\prec$ to denote the order relation of the latter.
We construct a $P$-tableau~$T$ as follows.
Consider the $\ell$ sinks of~$\Frak o$, i.e., the $\ell$
minimal elements of~$\Frak o$.
These are mutually non-adjacent in~$G$
and therefore they form a chain in~$P$.
Let the $\prec$-minimal element of this chain be the $(1,1)$ cell of~$T$
(i.e., the cell in the first row and the first column of~$T$).
Choose $k-1$ out of the remaining $\ell-1$ sinks and
arrange these in $\prec$-ascending order down the first column of~$T$.
Arrange the remaining elements along the first row of~$T$ as follows.
At each stage, the remaining elements form an induced subposet of~$\Frak o$.
The minimal elements of this subposet form a chain in~$P$.
Choose the $\prec$-minimal element of this chain to be the next
element in the first row of~$T$, and repeat this process
until all elements have been placed.

We must check that $T$ is a $P$-tableau.
We need only check that if $u$ and $v$ are consecutive elements
in the first row, then $u\not\succ v$.
If $u\succ v$, then in particular $u$ and~$v$ are $\prec$-comparable
and therefore non-adjacent in~$G$.
There are two cases.

\noindent {\it Case 1:
$u$ is not the element in the $(1,1)$ cell of~$T$.}
Then $u$ is a $<$-minimal element of some induced subposet~$\Frak o'$
of~$\Frak o$ and $v$ is a $<$-minimal element of $\Frak o' \setminus u$.
Now $<$-minimal elements of~$\Frak o'\setminus u$
are either $<$-minimal elements of~$\Frak o'$ or
else elements that cover~$u$ (in the $<$-ordering).
However, $v$ cannot cover~$u$ because $v$ and~$u$ are
non-adjacent in~$G$.  Hence $v$ is a $<$-minimal element of~$\Frak o'$,
and since $u$ is $\prec$-minimal among all $<$-minimal elements
of~$\Frak o'$, we must have $u\prec v$, a contradiction.

\noindent {\it Case 2: $u$ is in the $(1,1)$ cell of~$T$.}
Then $v$ is a $<$-minimal element of the subposet $\Frak o'$
obtained by deleting all elements in the first column of~$T$ from~$\Frak o$.
Again $v$ is either a $<$-minimal element of~$\Frak o$
or else $v$ covers one of the elements in the first column of~$T$.
But as in Case~1, $v$ cannot cover~$u$,
and if $u'$ is some other element in the first column of~$T$,
then $v \prec u \prec u'$ so $v$ cannot cover~$u'$ either.
Hence, arguing as before, $v$ is a $<$-minimal element of~$\Frak o$
and $u \prec v$, a contradiction.

Thus $T$ is indeed a $P$-tableau.
It is easy to see that $T$ induces~$\Frak o$.
The $\ell - 1 \choose k-1$ $P$-tableaux produced by
the procedure described above are all distinct
because they have distinct first columns.
It remains only to show that no other
$P$-tableau whose shape is a hook with $k$ elements
in the first column can induce~$\Frak o$.
It is clear that the first column
of any $P$-tableau~$T$ inducing~$\Frak o$
must consist
of $k$ $<$-minimal elements of~$\Frak o$ in $\prec$-ascending order.
We claim that the only possible element that can go in the $(1,1)$
cell of~$T$ is the
$\prec$-minimal element~$u$ of the $<$-minimal elements of~$\Frak o$.
For suppose that the element in the $(1,1)$ cell is some $v\ne u$.
Since $v$ is necessarily a $<$-minimal element of~$\Frak o$,
we have $u\prec v$, and therefore $u$ cannot be in the $(1,2)$
cell of~$T$.
We claim that $u$ cannot actually be anywhere in the first row
without violating the $P$-tableau condition.
Suppose that $w$ is the element in the $(1,2)$ cell.
Then $w$ is a $<$-minimal element of the subposet~$\Frak o'$
consisting of the elements in the first row of~$T$ excluding~$v$.
So either $w$ is a $<$-minimal element of~$\Frak o$
or $w$ covers some element~$w'$ in the first column of~$T$.
In the former case, $u\prec w$ by definition of~$u$,
and in the latter case, $u$ and~$w$ are $\prec$-comparable
(since $u\in \Frak o'$ and $w$ is $<$-minimal in~$\Frak o'$)
and we cannot have $w\prec u$ for then $w \prec u \prec u'$
would be a contradiction.
Either way, $u\prec w$ so $u$ cannot be in the $(1,3)$ of~$T$
either.  This argument can be continued inductively to show
that $u$ cannot be anywhere, a contradiction.
Following a similar argument,
we can show that the only possible way the elements in the first
row can be arranged is according to the algorithm given previously.
This completes the proof.\qed

\noindent{\it Proof of \acyclic.\enspace}
Define
$$c^G_\ell \defeq \suml a^G_\lambda.$$
By \mainthm,
$$\eqalignno{c^G_\ell &= \sum_T \sgn T,&\eqn\cglone\cr}$$
where the sum is over all special rim hook $P$-tableaux~$T$
with $\ell$ rim hooks.
We can break up the sum \cglone:
$$\eqalignno{c^G_\ell &= \sum_\mu \sum_T \sgn T,&\eqn\cgl\cr}$$
where the outer sum is over all shapes $\mu$
and the inner sum is over all special rim hook $P$-tableaux $T=(T',T'')$
having shape~$\mu$ and $\ell$ rim hooks.
We now claim that the inner sum in~\cgl\ vanishes unless $\mu$ is a hook.

To prove this claim, assume that $\mu$ is not a hook,
so that $\mu$ contains the cell $(2,2)$.
For any special rim hook $P$-tableau $T=(T',T'')$
with shape~$\mu$,
let $H_1(T)$ denote the (special) rim hook of~$T''$ containing $(1,1)$ and
let $H_2(T)$ denote the (special) rim hook of~$T''$ containing $(2,2)$.
It is easy to see that $H_1(T)\ne H_2(T)$.
Now fix any $P$-tableau~$T'$ of shape~$\mu$
and let $S$ be the set of all special rim hook $P$-tableaux
having shape~$\mu$ and $\ell$ rim hooks
and having the form $(T',T'')$ for some~$T''$.
We now define a sign-reversing involution~$\sigma$ on~$S$.
This will prove the claim.

Note that if $H_2(T)$ contains some cells in the first row,
say the cells $(1,m)$ through $(1,n)$ for some $m\le n$,
then all the cells $(1,1)$ through $(1,m-1)$ must belong
to~$H_1(T)$.  Note also that if $H_2(T)$ does not contain
any cells in the first row and if the rightmost cell of~$H_2(T)$
is $(2,m)$ for some~$m$, then all the cells
$(1,1)$ through $(1,m)$ must belong to~$H_1(T)$.
With this in mind we define $\sigma(T)$ as follows:
$\sigma(T)$ is exactly the same as~$T$ except
that $H_1(\sigma(T)) \ne H_1(T)$ and $H_2(\sigma(T))\ne H_2(T)$.
If $H_2(T)$ contains cells in the first row,
then transfer these first-row cells to~$H_1$, i.e., let $H_2(\sigma(T))$
equal $H_2(T)$ with the first-row cells of~$H_2(T)$ deleted,
and let $H_1(\sigma(T))$ equal $H_1(T)$ plus the first-row cells
of~$H_2(T)$.
If $H_2(T)$ does not contain any cells in the first row
and its rightmost cell is $(2,m)$,
then let $H_2(\sigma(T))$ equal $H_2(T)$
plus the first-row cells of~$H_1(T)$ to the right of $(1,m)$
(including $(1,m)$ itself)
and let $H_1(\sigma(T))$ equal $H_1(T)$ with all the cells
to the right of $(1,m)$ (including $(1,m)$ itself) deleted.
It is easy to check that the definition of $\sigma$ makes sense
and that it is an involution on the set~$S$.
It is sign-reversing because the sign of~$H_2$ is changed
(its height changes by one) but the signs of all the
other rim hooks remain unchanged.

Thus in \cgl\ we may restrict the outer sum to hook-shapes~$\mu$.
Now it is easy to see that for a given hook~$\mu$,
all special rim-hook $P$-tableaux having shape~$\mu$ and $\ell$ rim hooks
have the same sign, namely $(-1)^{k(\mu)-\ell}$,
where $k(\mu)$ is the number of cells in the first column of~$\mu$.
Moreover, we claim that the number of special rim hook tabloids
having shape~$\mu$ and $\ell$ rim hooks is $k(\mu)-1\choose \ell-1$.
For, because of the special condition that all rim hooks must
contain at least one cell in the first column,
the cells in the first row of~$\mu$ must all belong to the rim hook
containing $(1,1)$, and therefore
the set of special rim hook tabloids
is in bijection with the set of compositions of~$k(\mu)$ into $\ell$ parts
(just look at the way the special rim hooks subdivide the
set of cells in the first column).

Using the notation of \ptableau, we may use what we know
to rewrite \cgl:
$$\eqalignno{c^G_\ell &= \sum_\mu \sum_T \sgn T \cr
   &= \sum_\mu \sum_T (-1)^{k(\mu)-\ell}\cr
   &= \sum_\mu (-1)^{k(\mu)-\ell} {k(\mu)-1 \choose \ell-1} \pi_\ell \cr
   &= \sum_k (-1)^{k-\ell} {k-1 \choose \ell-1} \pi_\ell\cr
   &= \sum_k (-1)^{k-\ell} {k-1 \choose \ell-1}
                   \sum_m {m-1 \choose k-1} \kappa_m \cr
   &= \sum_m \kappa_m \sum_k (-1)^{k-\ell} {m-1 \choose k-1}
                                           {k-1 \choose \ell-1} \cr
   &= \sum_m \kappa_m {m-1 \choose \ell - 1}
                \sum_k (-1)^{k-\ell} {m-\ell \choose k-\ell} \cr
   &= \sum_m \kappa_m {m-1 \choose \ell - 1}
                \sum_k (-1)^k {m-\ell \choose k} \cr
   &= \sum_m \kappa_m {m-1 \choose \ell - 1} \delta_{m,\ell} \cr
   &= \kappa_\ell.\cr}$$
This proves the proposition.\qed

\sec Concluding Remarks

Optimistically, \mainthm\ could be
used to provide a combinatorial proof of
the Poset Chain Conjecture.
However, there is a caveat.
In the case when $P$ is an ordinal sum of antichains
of sizes $\nu_1, \nu_2, \ldots, \nu_m$,
one can show that $f^G_\mu = \nu_1! \nu_2! \cdots \nu_m!K_{\mu,\nu}$,
where the $K_{\mu,\nu}$ are the Kostka numbers.
The Poset Chain Conjecture in this case amounts to
the assertion that
if $K$ is the Kostka matrix,
then $K^{-1}K$ is a matrix with nonnegative coefficients.
While this is a trivial fact algebraically,
E\~gecio\~glu and Remmel state that
it is an open problem to prove this bijectively
using their combinatorial interpretation.
It seems therefore that \mainthm\ needs to be
supplemented by algebraic arguments for it to
be an effective tool for attacking
the full Poset Chain Conjecture.

\sec Acknowledgments

The author was supported in part by a National Science Foundation
postdoctoral fellowship and did part of the work for this paper
while a general member of the Mathematical Sciences Research
Institute.  Thanks also to Jaejin Lee, who spotted a minor error
in an earlier version of this manuscript.

\sec References

\item{1.} \"O. E\~gecio\~glu and J. B. Remmel,
A combinatorial interpretation of the inverse Kostka matrix,
\itc{Lin.\ Multilin.\ Alg.}\ {\bf 26} (1990), 59--84.

\item{2.} V. Gasharov,
Incomparability graphs of \tpof\ posets are $s$-positive.
\itc{Discrete Math.}\ {\bf 157} (1996), 193--197.

\item{3.} I. G. Macdonald,
``Symmetric Functions and Hall Polynomials,'' 2nd ed.,
Oxford University Press, New York, 1995.

\item{4.} B. Sagan,
``The Symmetric Group: Representations, Combinatorial  
Algorithms, and Symmetric Functions,''
Wadsworth \& Brooks/Cole, Pacific Grove, 1991.

\item{5.} R. P. Stanley, A symmetric function generalization of the
chromatic polynomial of a graph, \itc{Advances in Math.} {\bf 111} (1995),
166--194.

\item{6.} R. P. Stanley and J. Stembridge,
On immanants of Jacobi-Trudi
matrices and permutations with restricted position, \itc{J.~Combin.\
Theory (A)} {\bf 62} (1993), 261--279.

\bye